\documentclass[11pt]{article}

\usepackage{amssymb,latexsym}
\usepackage{amsmath}
\usepackage{citesort}

\title{On Positive Deformations of $^*$-Algebras\thanks{Talk given by
    Stefan Waldmann at the Mosh{\`{e}} Flato Conference 1999 in Dijon.}}

\author{{\bf 
        Henrique Bursztyn\thanks{henrique@math.berkeley.edu}
        } \\[0.5cm]
        Department of Mathematics\\
        UC Berkeley\\
        94720 Berkeley, CA, USA
        \\[1cm]
        {\bf 
        Stefan Waldmann\thanks{Stefan.Waldmann@ulb.ac.be}
        } \\[0.5cm]
        D{\'{e}}partement de Math{\'{e}}matique \\
        Campus Plaine, C. P. 213 \\
        Boulevard du Triomphe \\
        B-1050 Bruxelles \\
        Belgique 
        }

\date{October 1999}

\newcommand{\im} {{\mathrm i}}
\newcommand{\eu} {{\mathrm e}}

\newcommand{\cc} [1] {\overline {{#1}}}
\newcommand{\supp} {\mathop{{\mathrm {supp}}}}

\newcommand{\id} {{\mathsf {id}}}
\newcommand{\tr} {{\mathsf {tr}}}

\newcommand{\Hom} {{\mathsf {Hom}}}

\newcommand{\ring} [1] {{\mathsf {{#1}}}}

\newcommand{\SP} [1] {{\left\langle {{#1}} \right\rangle}}

\newcommand{\Unit} {\mathbf {1}}

\newcommand{\weyl}{\mathop{*_{\mbox{\tiny Weyl}}}}
\newcommand{\wick}{\mathop{*_{\mbox{\tiny Wick}}}}

\newtheorem{lemma} {Lemma} [section]
\newtheorem{proposition} [lemma] {Proposition}

\newtheorem{corollary} [lemma] {Corollary}
\newtheorem{definition}[lemma] {Definition}

\newenvironment{proof}{\small{\sc Proof:}}{{\hspace*{\fill} $\square$\\}}

\numberwithin{equation}{section}

\begin{document}

\maketitle

\begin{center}
Dedicated to the memory of Mosh{\'e} Flato
\end{center}

\vspace{0.5cm}

\begin{abstract}
Motivated by deformation quantization we consider $^*$-algebras over
ordered rings and their deformations: we investigate formal
associative deformations compatible with the $^*$-involution and
discuss a cohomological description in terms of a Hermitian Hochschild
cohomology. As an ordered ring allows for a meaningful definition of
positive functionals and as the formal power series with coefficients
in an ordered ring are again an ordered ring we define a deformation
to be positive if any positive linear functional of the undeformed
algebra can be deformed into a positive linear functional of the
deformed algebra. We discuss various examples and prove in particular
that star products on symplectic manifolds are positive deformations.
\end{abstract}

%
%

\section{Introduction}

In this paper we discuss several aspects of deformations of algebras
with a $^*$-involution which are defined over an ordered ring. Our
main example is deformation quantization as introduced by Bayen,
Flato, Fr{\o}nsdal, Lichnerowicz, and Sternheimer in \cite{BFFLS78},
see e.g.\  \cite{Wei94,Ste98} for recent surveys. In deformation
quantization the algebra of classical observables, the
smooth complex-valued functions $C^\infty (M)$ on a Poisson manifold,
is deformed into the quantum mechanical observable algebra by
introducing a new, $\hbar$-dependent associative product, the star
product, such that in zeroth order of $\hbar$ the star product
coincides with the pointwise product and the commutator yields in
first order $\im$ times the Poisson bracket. In deformation 
quantization the star product usually is regarded as formal,
i.e. based on formal power series in $\hbar$, and thus a star product
is a formal associative deformation in the sense of Gerstenhaber
\cite{GS88}. The existence of star products was shown by DeWilde and
Lecomte \cite{DL83b}, Fedosov \cite{Fed86,Fed94a}, and Omori, Maeda,
and Yoshioka \cite{OMY91} for the case of symplectic manifolds, and
by Kontsevich \cite{Kon97b} for the general case of a Poisson
manifold. The classification up to equivalence is due to Nest and
Tsygan \cite{NT95a,NT95b}, Bertelson, Cahen, and Gutt \cite{BCG97},
Deligne \cite{Del95}, Weinstein and Xu \cite{WX98}, and Kontsevich
\cite{Kon97b}. On the other hand the formal character of these star
products is certainly not sufficient from the physical point of view
whence one also has to study convergence properties with respect to
$\hbar$. This can be done either by starting with a formal star
product and then finding suitable subalgebras where one can make sense 
out of convergence, see
e.g.\  \cite{CGR90,CGR93,CGR94,CGR95,BBEW96b,BNW98a,BNW99a,BNPW98} and
references therein, or directly by formulating the deformation within
the framework of e.g.\  $C^*$-algebras from the beginning, which leads
e.g.~to Rieffel's notion of strict deformation quantization, see
\cite{Rie89} and Landsman's book \cite{Land98} for related approaches
and references. In many cases the formal approach can be seen as an
`asymptotic expansion' in $\hbar$ of a convergent situation, which can 
be made more precise in several examples as e.g. for cotangent
bundles \cite{BNW98a,BNW99a,BNPW98}, and it is very useful to 
think of formal deformations in that way even if it is not clear if
there exists a corresponding strict deformation. This leads to the
following programme, namely trying to find `formal analogs' of various 
techniques known in $C^*$-algebra theory in such a way that they can be 
viewed as `asymptotic expansions' and thus to understand the classical 
and semi-classical limit of these constructions which possibly opens
also an asymptotic point of view to Connes' non-commutative geometry
\cite{Con94}. This approach has
turned out to be interesting both from the physical and mathematical
point of view: in \cite{BW98a,BW97b} Bordemann and Waldmann have
introduced a notion of formal GNS representations for star products on 
pre-Hilbert spaces over the ring of formal power series by noticing
that the ring of real formal power series is in fact ordered and thus
allows to define positive functionals etc.\  as in the $C^*$-algebra
theory. This approach (re-)produces many physically interesting
results for the star product algebras in a purely formal framework, as 
e.g.\  for cotangent bundles the Schr{\"o}dinger-like differential
operator representations, Aharonov-Bohm-like effects and
representations on sections of line bundles when magnetic monopoles
are present \cite{BNW98a,BNW99a,BNPW98}. Also thermodynamical KMS
states and their representations have been investigated, see
\cite{BFLS84,BL85} as well as \cite{BRW98a,BRW98b} and \cite{Wal99a}
for a related `formal' Tomita-Takesaki theory. On the other hand the
framework of $^*$-algebras and pre-Hilbert spaces over ordered rings
seems to be mathematically interesting and turns out to have a rich
structure. In \cite{BuWa99a} we started to develop the notions of
algebraic Rieffel 
induction and formal Morita equivalence for such $^*$-algebras in
order to understand the (semi-)classical limits of the corresponding
constructions for $C^*$-algebras, as e.g.\  used by Landsman
\cite{Land98} to describe quantum analogous of phase space
reduction. As this also has been done in the framework of deformation
quantization, see e.g.\  \cite{Fed98a,BHW99a}, it would be very
interesting to compare these approaches. During our investigations we
needed the notion of a \emph{positive deformation} of a $^*$-algebra
in the sense that all positive linear functionals can be deformed 
into positive linear functionals of  the deformed algebra. This is in
some sense `dual' to Landsman's notion of a positive quantization map,
but less restrictive, and it seems to be more suitable for our
purposes.

In this paper we shall discuss the deformations of $^*$-algebras over
ordered rings in general and apply our results to the main examples of 
deformation quantization obtaining the following results: after
reviewing briefly the notions of ordered rings and $^*$-algebras  and
discussing the notion of `sufficiently many positive linear
functionals' in Section~\ref{OrderedSec}, we set-up the deformation
programme for $^*$-algebras in Section~\ref{DefSec} where the `real
part' of the Hochschild cohomology is relevant for Hermitian
deformations, i.e.\  those deformations which do not deform the
$^*$-involution. In Section~\ref{PosDefSec} we motivate the notion of
a positive deformation and discuss some basic constructions proving in
particular that positive deformations of $^*$-algebras with
sufficiently many positive linear functionals still have sufficiently
many positive linear functionals. Finally, in Section~\ref{QuantSec}
we prove that star products on symplectic manifolds are positive
deformations and that on arbitrary Poisson manifolds star products
have sufficiently many positive linear functionals.

%
%

\section{$^*$-Algebras over ordered rings}
\label{OrderedSec}

In this section we shall recall some basic definitions and examples
concerning ordered rings and $^*$-algebras. For further references
and details see e.g.\  \cite{BW98a,BuWa99a}.
Let $\ring R$ be an \emph{ordered ring}, i.e.\  an associative,
commutative and unital ring with a subset $\ring P \subset \ring R$
such that $\ring R$ is the disjoint union 
$\ring R = - \ring P \cup \{0\} \cup \ring P$ and 
$\ring P \cdot \ring P \subseteq \ring P$,
$\ring P + \ring P \subseteq \ring P$. The elements in $\ring P$ are
called positive. Clearly $\ring R$ has characteristic zero and no
zero divisors. Then consider the quadratic ring extension 
$\ring C = \ring R(\im)$ where $\im^2 := -1$. Complex conjugation 
$z = a + \im b \mapsto \cc z :=a - \im b$ with $a,b \in \ring R$ is
defined as usual and clearly $z \in \ring R \subset \ring C$ iff 
$z = \cc z$ and $\cc z z \ge 0$.

Let $\mathcal A$ be an associative algebra over $\ring C$ endowed with 
a $^*$-involution, i.e. an involutive $\ring C$-antilinear
anti-automorphism $^*: \mathcal A \to \mathcal A$, then $\mathcal A$
is called a \emph{$^*$-algebra} over $\ring C$. As usual we define
\emph{Hermitian}, \emph{unitary} and \emph{normal elements} in 
$\ring C$. Given such a $^*$-algebra we  
define a linear functional $\omega: \mathcal A \to \ring C$ to be
\emph{positive} if $\omega(A^*A) \ge 0$ for all $A \in \mathcal A$.
If $\mathcal A$ is unital then a positive linear functional is called
a state if $\omega(\Unit) = 1$. For any positive linear functional one 
has the Cauchy-Schwarz inequality 
$\omega(A^*B)\cc{\omega(A^*B)} \le \omega(A^*A)\omega(B^*B)$ 
and $\omega (A^*B) = \cc{\omega (B^*A)}$ for 
$A, B \in \mathcal A$. Next we consider positive algebra elements in
$\mathcal A$. There are essentially two ways to define positivity 
\cite[Def.~2.3]{BuWa99a}:
firstly, we call $A \in \mathcal A$ \emph{algebraically positive} if 
$A = a_1 A_1^*A_1 + \cdots + a_n A_n^*A_n$, where $n \in \mathbb N$,
$a_i \ge 0$, and $A_i \in \mathcal A$ for $i = 1, \ldots, n$, and we
denote the set of all algebraically positive elements by 
$\mathcal A^{++}$. Secondly, we define a Hermitian element $A$ to be
\emph{positive} if for all positive linear functionals 
$\omega: \mathcal A \to \ring C$ one 
has $\omega(A) \ge 0$, and we denote the positive elements by 
$\mathcal A^+$. Then clearly 
$\mathcal A^{++} \subseteq \mathcal A^+$ and both sets are convex
cones with the property that for any $C \in \mathcal A$ one has
$C^* \mathcal A^{++} C \subseteq \mathcal A^{++}$ as well as
$C^* \mathcal A^+ C \subseteq \mathcal A^+$.

The basic examples of such $^*$-algebras are given by $\ring C$
itself, in which case $\ring C^+ = \ring C^{++}$ coincides with the
non-negative elements in $\ring R$. More generally, one can consider 
$M_n (\ring C)$, the $n \times n$ matrices over $\ring C$ with the
usual $^*$-structure. As  
in the case of $\ring C = \mathbb C$ one can show that the positive
linear functionals can be written as $A \mapsto \tr(\varrho A)$ where
$\tr$ is the trace functional and $\varrho$ is a Hermitian matrix such 
that $\SP{v, \varrho v} \ge 0$ for all $v \in \ring C^n$, where
$\SP{\cdot,\cdot}$ denotes the usual Hermitian product on 
$\ring C^n$. Moreover, $A = A^* \in M_n (\ring C)^+$ if  and only if
$\SP{v, Av} \ge 0$ for all $v \in \ring C^n$. In the case of
deformation quantization we consider the space of complex-valued
smooth functions $C^\infty_0 (M)$ with compact support (only
for technical reasons) on a Poisson manifold as the space of classical 
observables with the obvious $^*$-algebra structure. Then the positive 
linear functionals are positive Borel measures with finite volume for
every compact set and $C^\infty_0 (M)^+$ consists of functions with
$f(x) \ge 0$ for all $x \in M$ as expected. Nevertheless we should
also mention that the above definitions lead to a somehow
`pathological' characterization of positive elements in the case of
polynomial algebras, see e.g.\  \cite{Pow71}. The final example we
would like to mention is the Grassmann algebra 
$\bigwedge (\ring C^n)$, where we define a $^*$-involution by requiring
$\Unit^* =\Unit$ and $e_i^* = e_i$, where $e_1, \ldots, e_n$ is the
canonical basis of $\ring C^n$. Then one easily shows that 
$\omega: \bigwedge (\ring C^n) \to \ring C$ is positive if and only if 
$\omega (\Unit) \ge 0$ and 
$\omega (e_{i_1} \wedge \cdots \wedge e_{i_r}) = 0$ for all $r \ge 1$
and $i_1, \ldots, i_r$. Thus in this example one has (up to
normalization) only one positive linear functional.

To avoid such examples we state the following definition
\cite[Def.~2.7]{BuWa99a}: a 
$^*$-algebra $\mathcal A$ has sufficiently many positive linear
functionals if for every Hermitian element $H \ne 0$ there exists a
positive linear functional $\omega$ such that $\omega (H) \ne 0$. 
Clearly $M_n(\ring C)$ and $C^\infty_0 (M)$ have sufficiently many
positive linear functionals but $\bigwedge (\ring C^n)$ does not.
This condition turns out to be rather strong as it implies already
some very $C^*$-algebra like properties: one can show that if
$\mathcal A$ has sufficiently many positive linear functionals then,
e.g. $A^*A = 0$ implies $A = 0$ as well as $H^n = 0$ for a normal
element $H$ implies $H = 0$. Moreover, there exists a faithful
$^*$-representation on a pre-Hilbert space over $\ring C$ which can be 
constructed by means of the GNS construction, see \cite{BuWa99a}.

%
%

\section{Deformations of $^*$-algebras}
\label{DefSec}

Before studying deformations of $^*$-algebras over ordered rings one
notices that if $\ring R$ is an ordered ring then the ring of formal
power series $\ring R[[\lambda]]$ is again ordered in a canonical way, 
namely one defines $a = \sum_{r=r_0}^\infty \lambda^r a_r$ to be 
positive if $a_{r_0} > 0$. Thus the deformation programme stays in the 
same category of $^*$-algebras over ordered rings (and their quadratic 
ring extension by `$\im$'). In the case of deformation quantization
the deformation parameter $\lambda$ plays the role of $\hbar$ and may
be substituted $\lambda \mapsto \hbar$ in \emph{convergent}
situations.

Denote by $\mu_0: \mathcal A \otimes \mathcal A \to \mathcal A$ the
$\ring C$-bilinear associative product of $\mathcal A$ and by 
$I_0: \mathcal A \to \mathcal A$ the $^*$-involution, which we shall
frequently refer to as the `classical' structures. Then a 
\emph{formal $^*$-algebra deformation} of $\mathcal A$ is a formal
deformation $\mu = \mu_0 + \sum_{r=1}^\infty \lambda^r \mu_r$ of
$\mu_0$ to an associative $\ring C[[\lambda]]$-bilinear product 
$\mu$ for $\mathcal A[[\lambda]]$ together with a formal deformation
$I = I_0 + \sum_{r=1}^\infty \lambda^r I_r$ of the $^*$-involution
$I_0$ to a $\ring C[[\lambda]]$-antilinear involutive
anti-automorphism of $\mu$. Two such deformations $(\mu, I)$ and
$(\tilde\mu, \tilde I)$ are called \emph{equivalent} if there exists a
formal series $T = \id + \sum_{r=1}^\infty \lambda^r T_r$ of 
$\ring C$-linear maps $T_r: \mathcal A \to \mathcal A$ such that 
$T(\mu (A \otimes B)) = \tilde\mu(TA \otimes TB)$ and 
$TI = \tilde IT$. Due to the physical interpretation in deformation
quantization of the Hermitian elements with respect to the classical
$^*$-involution as observables (which are the real-valued functions),
we are mainly interested in those $^*$-algebra deformations where
the classical $^*$-involution is \emph{not} deformed, since the
classical observables should be observable in the `quantized' theory,
too. This particular kind of $^*$-algebra deformations shall be called
\emph{Hermitian deformations} of $\mathcal A$ (or symmetric, see
\cite{BFFLS78} for this definition in the context of star
products). Equivalence of Hermitian deformations now simply means that
there is a real equivalence transformation.

We shall now briefly sketch the deformation theory for Hermitian
deformations in terms of a Hermitian Hochschild cohomology for
$\mathcal A$, see \cite{GS88} for the general cohomological approach
to deformations of algebras. On the Hochschild complex 
$C^\bullet (\mathcal A) = 
\bigoplus_{r=0}^\infty \Hom(\mathcal A^{\otimes n}, \mathcal A)$
one defines an action of the $^*$-involution 
(we sometimes use again the symbol $^*$ instead of $I_0$) by
\begin{equation}
\label{InvHochschild}
    \varphi^* (A_1 \otimes \cdots \otimes A_n) 
    := (\varphi (A_n^* \otimes \cdots \otimes A_1^*))^*,
\end{equation}
where $A_1, \ldots, A_n \in \mathcal A$ and $\varphi$ is a
$n$-cochain. Clearly $\varphi^*$ is again a $n$-cochain and the map
$\varphi \mapsto \varphi^*$ is $\ring C$-antilinear and involutive.
Note that we have also reverted the order of the arguments. Then one
finds by a straightforward computation the following relation for the
Gerstenhaber product $\varphi \diamond \psi$ of two cochains 
$\varphi \in C^n(\mathcal A)$, $\psi \in C^m(\mathcal A)$:
\begin{equation}
\label{InvGerstenhaber}
    (\varphi \diamond \psi)^* 
    = (-1)^{(n-1)(m-1)} \; \varphi^* \diamond \psi^*.
\end{equation}
The fact that $^*$ is an anti-automorphism of $\mu_0$ is expressed by
the relation
\begin{equation}
\label{InvMu}
    \mu_0^* = \mu_0,
\end{equation}
whence we obtain for the Hochschild differential with respect to
$\mu_0$ for a $n$-cochain $\varphi$
\begin{equation}
\label{InvDelta}
    (\delta \varphi)^* = (-1)^{n-1} \delta \varphi^*,
\end{equation}
where $\delta$ is the Gerstenhaber bracket with $\mu_0$, i.e. the
graded $\diamond$-commutator 
$\delta\varphi = (-1)^{n-1} [\mu_0, \varphi]$.  
We call a cochain $\varphi$ \emph{Hermitian} 
if $\varphi^* = \varphi$ and anti-Hermitian if 
$\varphi^* = - \varphi$, and denote the \emph{Hermitian cochains} by 
$C^\bullet_H (\mathcal A)$. They are a $\mathbb Z$-graded 
$\ring R$-submodule of $C^\bullet (\mathcal A)$. For a $0$-cochain
this coincides with the notion of a Hermitian algebra element while a
$1$-cochain is Hermitian if $\varphi (A^*) = \varphi (A)^*$ for all 
$A \in \mathcal A$. Moreover, for $\varphi \in C^n_H (\mathcal A)$ the 
cochain $\delta \varphi$ is Hermitian if $n$ is odd and anti-Hermitian 
if $n$ is even. This observation leads to the following definition of
the Hermitian Hochschild cohomology of $\mathcal A$. Define
$Z^n_H(\mathcal A)$ to be the set of Hermitian cocycles
and define $B^n_H(\mathcal A)$ to be the set of Hermitian 
coboundaries $\varphi = \delta \psi$ where $\psi$ is Hermitian if $n$
even and anti-Hermitian if $n$ is odd. Then 
$H^n_H(\mathcal A) := Z^n_H(\mathcal A) \big/ B^n_H(\mathcal A)$
denotes the $n$-th \emph{Hermitian Hochschild cohomology group}. As an 
ordered ring $\ring R$ has characteristic zero and no zero-divisors it 
is reasonable to assume in addition $\frac{1}{2} \in \ring R$ since
one may pass to an appropriate quotient anyway. Now if 
$\frac{1}{2} \in \ring R$ then we can decompose any cochain 
$\ring R$-linearly into a Hermitian and anti-Hermitian part, i.e.
\begin{equation}
\label{HermitianCochain}
   \varphi = \frac{1}{2} (\varphi + \varphi^*)
             + \im \frac{1}{2\im} (\varphi - \varphi^*),
\end{equation}
and an easy check using (\ref{InvDelta}) shows that the same holds on
the level of cohomology whence we have the following $\ring R$-linear
isomorphism 
\begin{equation}
\label{HermitianCoho}
    \begin{array}{c}
    H^\bullet (\mathcal A) 
    \cong H^\bullet_H (\mathcal A) \oplus 
    \im H^\bullet_H (\mathcal A) 
    \\
    \left[\varphi\right] 
    \mapsto \left[\frac{1}{2}(\varphi + \varphi^*)\right]
    + \im \left[\frac{1}{2\im}(\varphi - \varphi^*)\right].
    \end{array}
\end{equation}

Now let us consider a deformation 
$\mu = \sum_{r=0}^\infty \lambda^r \mu_r$ of
$\mu_0$. Then the associativity condition for $\mu$ is well-known to be 
equivalent to $[\mu,\mu] = 0$ which yields $[\mu_0,\mu_0] = 0$ in
zeroth order and in higher orders of $\lambda$ one has
\begin{equation}
\label{AssObstr}
    \delta\mu_r = -\frac{1}{2}\sum_{s=1}^{r-1} [\mu_s,\mu_{r-s}].
\end{equation}
To obtain in addition a Hermitian deformation we need 
\begin{equation}
\label{HermObstr}
    \mu_r^* = \mu_r
\end{equation}
for all $r$ since $\lambda$ is `real' in $\ring C[[\lambda]]$. Thus we
end up with the following observation that if one can solve both
associativity and Hermiticity up to order $r-1$, then the obstructions
to solve (\ref{AssObstr}) and (\ref{HermObstr}) are given by 
$H^3_H (\mathcal A)$. If on the other hand one is able to solve the
cohomological equation (\ref{AssObstr}) for $\mu_r$ then there is even 
a Hermitian solution since an easy computation shows that the right
hand side in (\ref{AssObstr}) is 
anti-Hermitian and thus by (\ref{HermitianCoho}) one may simply take
the Hermitian part $\frac{1}{2}(\mu_r + \mu_r^*)$ instead of
$\mu_r$. Nevertheless note that this is in general not cohomologous to 
$\mu_r$ and thus there may be associative deformations of $\mu_0$
which are \emph{not} equivalent to a Hermitian deformation. This
happens e.g. for star products on symplectic manifolds with
non-trivial second deRham cohomology, see \cite{Neu99}. We summarize
these results in the following proposition:
\begin{proposition}
\label{HermHochProp}
Let $\mathcal A$ be a $^*$-algebra and let 
$\mu_0 + \lambda \mu_1 + \ldots + \lambda^{r-1}\mu_{r-1}$ be a
Hermitian deformation up to order $r-1$. Then the obstruction space
for a Hermitian deformation up to order $r$ is given by 
$H^3_H (\mathcal A)$. If there exists an associative deformation up to 
order $r$ then there also exists an Hermitian deformation up to order
$r$.
\end{proposition}

%
%

\section{Positive deformations}
\label{PosDefSec}

Let us now discuss the main object of this work, the definition of
positive deformations of $^*$-algebras. In order to motivate this
let us first remember the following example of 
\cite[Sect.~2, p.~555]{BW98a}: Consider the phase space 
$\mathbb R^{2n}$ with the canonical coordinates 
$q^1, \ldots, q^n, p_1, \ldots, p_n$ and the symplectic form
$\omega = \sum_k dq^k \wedge dp_k$. Then the Weyl-Moyal star product,
see \cite{BFFLS78}, is given by
\begin{equation}
\label{WeylMoyal}
    f \weyl g = \mu_0 \circ \eu^{\frac{\im\lambda}{2} 
                \sum_k (\partial_{q^k}\otimes \partial_{p_k} 
                - \partial_{p_k} \otimes \partial_{q^k})}
                f \otimes g,
\end{equation}
where $f, g \in C^\infty (\mathbb R^{2n})[[\lambda]]$, and $\mu_0$ is
the pointwise product. Clearly 
$\cc {(f \weyl g)} = \cc g \weyl \cc f$
and thus $\weyl$ is a Hermitian deformation of the pointwise product
with the complex conjugation as $^*$-involution. Consider the
Hamiltonian $H = \sum_k(q^kq^k + p_kp_k)$ of the harmonic oscillator
and let 
$\delta: 
C^\infty (\mathbb R^{2n})[[\lambda]] \to \mathbb C[[\lambda]]$
be the $\delta$-functional at $0$ which is classically a positive
linear functional. Then
\begin{equation}
\label{deltaNotpos}
    \delta(H\weyl H) = - \frac{\lambda^2}{2} < 0
\end{equation}
shows that $\delta$ is no longer a positive linear functional for the
Weyl-Moyal star product. On the other hand it is clear that the
classical limit of a $\ring C[[\lambda]]$-linear functional of a
deformed $^*$-algebra $\mathcal A[[\lambda]]$ has to be positive with
respect to the undeformed $^*$-algebra structure. This motivates the
following definition:
\begin{definition}
Let $\mathcal A$ be a $^*$-algebra over $\ring C$ and let 
$(\mathcal A[[\lambda]], \mu, I)$ be a $^*$-algebra deformation of
$\mathcal A$. 
\begin{enumerate}
\item The deformation is called positive if for any classically
      positive $\ring C$-linear functional $\omega_0$ there exist
      $\ring C$-linear maps $\omega_r: \mathcal A \to \ring C$, 
      $r\ge 1$, such that the $\ring C[[\lambda]]$-linear functional
      $\omega = \sum_{r=0}^\infty \lambda^r \omega_r: 
      \mathcal A[[\lambda]] \to \ring C[[\lambda]]$ is positive with
      respect to the deformed $^*$-algebra structure.
\item The deformation is called strongly positive if any classically
      positive linear functional is positive for the deformed
      $^*$-algebra, too.
\end{enumerate}
\end{definition}
Clearly positivity of a $^*$-algebra deformation is a property of the
whole corresponding equivalence class of $^*$-algebra deformations as
one can pull-back positive functionals with the equivalence
transformation which yields again a positive functional without
changing the classical limit. Nevertheless, this is no longer true for 
strongly positive deformations in general and it remains to find an
appropriate notion of equivalence here.

As a first application we observe that having sufficiently many
positive linear functionals is a property which is preserved under
positive deformations:
\begin{proposition}
\label{ManyPosProp}
Let $\mathcal A$ be a $^*$-algebra over $\ring C$ having sufficiently
many positive linear functionals. Then every positive deformation of
$\mathcal A$ has sufficiently many positive linear functionals. 
\end{proposition}
\begin{proof}
Let $(\mathcal A[[\lambda]], \mu, I)$ be a positive deformation of
$(\mathcal A, \mu_0, I_0)$ and let 
$A = \sum_{r=r_0}^\infty \lambda^r A_r$ be a non-zero Hermitian
element with respect to $I$ with first non-vanishing order $r_0$. Then 
$A_{r_0} = I_0(A_{r_0})$ is Hermitian with respect to the classical
$^*$-involution and thus we find a positive $\ring C$-linear
functional $\omega_0$ of $\mathcal A$ with 
$\omega_0 (A_{r_0}) \ne 0$. Moreover, we find 
$\omega_r: \mathcal A \to \ring C$ such that 
$\omega = \sum_{r=0}^\infty \lambda^r \omega_r$ is a positive
$\ring C[[\lambda]]$-linear functional of $\mathcal A[[\lambda]]$ with 
respect to the deformed $^*$-algebra structure. But then clearly
$\omega(A) \ne 0$, proving the proposition.
\end{proof}

As $^*$-algebras with sufficiently many positive linear functionals
behave in many aspects almost like $C^*$-algebras we are mainly
interested in such positive deformations.

The main difficulty in discussing the existence of (strongly) positive 
deformations is that one has to deal with inequalities instead of
equalities which seems to exclude a suitable cohomological
approach. It may even be difficult to decide whether a given
$^*$-algebra deformation, as e.g. the Weyl-Moyal star product, is
positive or not. Nevertheless in case of a Hermitian deformation we
have the following simple but rather useful criterion:
\begin{lemma}
Let $(\mathcal A[[\lambda]], \mu, I_0)$ be a Hermitian deformation of
a $^*$-algebra $\mathcal A$ over $\ring C$. If for all 
$A \in \mathcal A$ and $r \ge 1$ one has 
$\mu_r(A^*\otimes A) \in \mathcal A^+$ then the deformation $\mu$ is
strongly positive.
\end{lemma}
\begin{proof}
Thanks to \cite[Lem.~A.5]{BNW99a} it is sufficient to check the
positivity of a $\ring C[[\lambda]]$-linear functional 
$\omega: \mathcal A[[\lambda]] \to \ring C[[\lambda]]$ on elements in
$\mathcal A$ alone. But if $\omega_0:\mathcal A \to \ring C$ is a
positive linear functional then clearly 
$\omega_0 (\mu(A^* \otimes A)) 
= \sum_{r=0}^\infty \lambda^r \omega_0 (\mu_r (A^* \otimes A)) \ge 0$.
\end{proof}

Based on this lemma we have the following class of strongly positive
deformations of $^*$-algebras constructed by using commuting
derivations:
\begin{lemma}
\label{LikeWickLem}
Let $\mathcal A$ be a $^*$-algebra over $\ring C = \ring R(\im)$, where 
we assume $\mathbb Q \subseteq \ring R$, and let $D_1, \ldots, D_n$ be 
derivations of $\mathcal A$ such that $[D_i, D_j] = 0 = [D_i, D_j^*]$
for all $i,j = 1, \ldots, n$. Then
\begin{equation}
\label{LikeWickProd}
    \mu = \mu_0 \circ \eu^{\lambda \sum_k D_k \otimes D_k^*}
\end{equation}
defines a strongly positive Hermitian deformation of $\mathcal A$. 
\end{lemma}
\begin{proof}
First we note that the $\exp$-series is well-defined thanks 
$\mathbb Q \subseteq \ring R$. As already observed , see
(\ref{InvDelta}), the $\ring C$-linear maps $D_i^*$ are again
derivations and thus all occuring derivations commute which implies
that (\ref{LikeWickProd}) is associative. Clearly
$\mu$ is Hermitian and one has 
$\mu_r(A^*\otimes A) \in \mathcal A^{++} \subseteq \mathcal A^+$
for all $r\ge 1$ whence we can apply the above lemma.
\end{proof}

%
%

\section{Positive deformations in deformation quantization}
\label{QuantSec}

Let us now discuss the situation in deformation quantization. The main 
example for a deformation of the type as in Lemma~\ref{LikeWickLem} is 
the Wick product $\wick$ for $\mathbb C^n$: denote by 
$z^1, \ldots, z^n$, $\cc z^1, \ldots, \cc z^n$ the global (anti-)
holomorphic coordinates then $\wick$ is defined for 
$f, g \in C^\infty (\mathbb C^n)[[\lambda]]$ by
\begin{equation}
\label{WickProd}
    f \wick g = \mu_0 \circ 
    \eu^{2\lambda \sum_k \partial_{z^k} \otimes \partial_{\cc z^k}} 
    f \otimes g,
\end{equation}
which is clearly of the above form and thus strongly positive and
Hermitian. See \cite{Kar96,Kar98a,Kar99,BW97a} for a more general
treatment of such star products with separation of variables on
K{\"a}hler manifolds. From the above remark on equivalent Hermitian
deformations and the well-known fact that on 
$\mathbb C^n \cong \mathbb R^{2n}$ all star products are
equivalent we conclude that the Weyl-Moyal star product $\weyl$ is an
example of a positive but \emph{not} strongly positive
deformation. More explicitly, we have the following: Using the explicit
(and real) equivalence transformation $T = \eu^{\lambda\Delta}$, where  
$\Delta = \sum_k \partial_{z^k}\partial_{\cc z^k}$ is the Laplacian
and $f \weyl g = T^{-1} (Tf \wick Tg)$, we conclude that for any
classically positive linear functional $\omega_0$ the 
$\mathbb C[[\lambda]]$-linear functional 
$\omega = \omega_0 \circ T = \omega_0 \circ \eu^{\lambda \Delta}$ is
positive with respect to $\weyl$. But note also that for some
particular $\omega_0$ there may be no need to deform it in order to
obtain positivity, see e.g.\  \cite{BNW98a,BNW99a}.

A local version of the above construction can even be used to show that 
all Hermitian star products on a symplectic manifold are positive
deformations:
\begin{proposition}
\label{StarProdPos}
Let $(M, *)$ be a symplectic manifold with a Hermitian star
product. Then $*$ is a positive deformation.
\end{proposition}
\begin{proof}
Let $\omega_0: C^\infty_0 (M) \to \mathbb C$ be a classically positive 
linear functional. Choose a locally finite open cover
$\{O_\alpha\}_{\alpha \in I}$ of $M$ by contractable charts and let
$\{\chi_\alpha\}_{\alpha \in I}$ be a subordinate `quadratic'
partition of unity, i.e. 
$\sum_\alpha \cc\chi_\alpha \chi_\alpha = 1$. Endow each $O_\alpha$
with a local Wick star product $*_\alpha$ and let $T_\alpha$ be a
local and real equivalence transformation between $*|_{O_\alpha}$ and
$*_\alpha$, i.e. $T_\alpha (f * g) = T_\alpha f *_\alpha T_\alpha g$
on $O_\alpha$. Such an equivalence transformation exists since
$O_\alpha$ is contractable, see e.g.\  \cite{BCG97}. Then it is 
easy to check that the $\mathbb C[[\lambda]]$-linear functional 
$\omega: C^\infty_0 (M)[[\lambda]] \to \mathbb C[[\lambda]]$ defined
by
\begin{equation}
\label{DeformedOmega}
    \omega (f) := \sum_\alpha \omega_0 \left( T_\alpha \left(
                  \cc \chi_\alpha * f * \chi_\alpha \right)\right)
\end{equation}
is well-defined, positive with respect to $*$, and a deformation of
$\omega_0$. 
\end{proof}

Since $C^\infty_0 (M)$ (or $C^\infty (M)$, respectively) has
sufficiently many positive linear functionals we obtain from
Prop.~\ref{ManyPosProp} and Prop.~\ref{StarProdPos} the following
useful corollary:
\begin{corollary}
Let $(M,*)$ be a symplectic manifold with a Hermitian star
product. Then the $\mathbb C[[\lambda]]$-algebra 
$(C^\infty_0 (M)[[\lambda]], *)$ has sufficiently many
positive $\mathbb C[[\lambda]]$-linear functionals.
\end{corollary}

The case of a Hermitian star product on an arbitrary Poisson manifold
seems to be more involved since the proof of Prop.~\ref{StarProdPos}
uses at several points that the Poisson structure is symplectic. Thus it
remains an open question whether or which star products on arbitrary
Poisson manifolds are positive deformations. Nevertheless we can show
directly the weaker property that a Hermitian star product on a
Poisson manifold has sufficiently many positive linear functional.
\begin{proposition}
\label{PoissonPos}
Let $(M, *)$ be a Poisson manifold with a Hermitian star product. Then 
the $\mathbb C[[\lambda]]$-algebra $(C^\infty_0 (M)[[\lambda]], *)$
has sufficiently many positive $\mathbb C[[\lambda]]$-linear
functionals. 
\end{proposition}
\begin{proof}
Let 
$0 \ne f = \sum_{r=0}^\infty \lambda^r f_r 
\in C^\infty_0 (M)[[\lambda]]$ 
be Hermitian then we may assume that $f_0 \ne 0$. Thus we find a point 
$x \in M$ and an open neighborhood $U \subseteq M$ of $x$ such that 
$f|_{U}$ is either strictly positive or strictly negative. Choose now
a smooth density 
$\mu \in \Gamma^\infty_0 (|\!\bigwedge^n\!|\, T^*M)$
with $\mu(x) > 0$, $\mu\ge 0$, and $\supp \mu \subset U$. Then it follows
by the same argument as in \cite[Lem.~2]{BW98a} that the functional 
$\omega: g \mapsto \int_M g\mu$ is positive with respect to any star
product on $M$ and clearly we have $\omega(f) \ne 0$.
\end{proof}

Let us conclude with a few remarks: It is clear that if a deformation
$\omega$ of a classically positive linear functional $\omega_0$ into a 
positive $\ring C[[\lambda]]$-linear functional of a deformed
algebra $\mathcal A[[\lambda]]$ exists, the higher orders $\omega_r$
of the deformation are 
not necessarily unique. If for example the functional $\omega_0$ is
faithful, i.e. if $\omega_0 (A^*A) > 0$ for all $A \ne 0$, then one can 
add \emph{any} real linear functionals $\omega_r$ in higher orders and 
does not loose positivity of $\omega$. This raises the question
whether in a particular case there are in some sense `natural' or
`minimal' corrections $\omega_r$ which deform $\omega_0$ into a
positive linear functional. In the case of star products one can
impose further conditions on $\omega_r$ by requiring e.g. continuity
with respect to the canonical locally convex topologies of smooth
functions. Moreover, one can demand that the corrections should not
increase the support of  
the classical functional $\omega_0$, see \cite{Wal99a} for an
extensive discussion on the involved locality structures. But even
then a tremendous variety remains, since one can pull-back deformed
positive linear functionals by $^*$-automorphisms of the deformed
algebra with the identity as classical limit. Then one ends up again
with a deformation of the same classically positive linear
functional. Hence a classification of positive deformations of a given 
classically positive linear functional modulo the above action of
$^*$-automorphisms would be highly desirable.

\section*{Acknowledgements}

We would like to thank Martin Bordemann, Nikolai Neumaier, and Alan
Weinstein for valuable discussions. One of us (S.W.) acknowleges
support of the Action de Recherche Concert{\'e}e de la Communaut{\'e}
fran{\c c}aise de Belgique.

\begin{small}

\end{small}
\end{document}